\documentclass[a4paper]{article}
\usepackage{graphicx}
\usepackage{onecolceurws}

\usepackage{hyperref}
\usepackage{amsmath}
\usepackage{amsthm}
\usepackage{amssymb}
\usepackage{pbox}

\usepackage{url}
\newcommand\thlink[1]{\href{http://us.metamath.org/mpegif/#1.html}{\textsf{#1}}}
\newtheorem{theorem}{Theorem}

\title{Formalization of the prime number theorem and Dirichlet's theorem}

\author{
Mario Carneiro \\ Pure and Applied Logic program \\
                Carnegie Mellon University, Pittsburgh PA, USA \\ di.gama@gmail.com
}

\institution{}

\begin{document}
\maketitle

\begin{abstract}
We present the formalization of Dirichlet's theorem on the infinitude of primes in arithmetic progressions, and Selberg's elementary proof of the prime number theorem, which asserts that the number $\pi(x)$ of primes less than $x$ is asymptotic to $x/\log x$, within the proof system Metamath.
\end{abstract}
\vskip 32pt

\section{Introduction}\label{sec:intro}
Dirichlet's theorem, or the Dirichlet prime number theorem, states that for any $N\in\mathbb N$ and $A\in\mathbb Z$ such that $\gcd(A,N)=1$, there are infinitely many primes in the progression $A+kN$, or equivalently there are infinitely many $p_k\equiv A\pmod N$. Euler was the first to make progress on this theorem, proving it in the case $A=1$, and it was shown in full generality by Dirichlet in 1837 \cite{dirichlet}.

The prime number theorem gives an overall order of growth of the number of primes less than $x$. Letting $\pi(x)$ denote the number of primes in the interval $[1,x]$ (where $x$ is not necessarily an integer), the prime number theorem asserts that $\pi(x)\sim\frac x{\log x}$. This was first conjectured by Legendre in 1797, and was first proven using complex analysis and the zeta function by Jacques Hadamard and Charles Jean de la Vall\'{e}e-Poussin in 1896. Two ``elementary'' proofs were later discovered by Selberg and Erd\H{o}s in 1949 \cite{selberg,erdos}.

The first formal proof of the prime number theorem was written by Jeremy Avigad et. al. in 2004 \cite{avigad}, in the Isabelle proof language, following Selberg's proof, a landmark result for mathematics formalization. Hadamard and Vall\'{e}e-Poussin's proof was formalized in HOL Light by John Harrison in 2009 \cite{harrison}. John Harrison also later formalized Dirichlet's theorem in HOL Light, in 2010 \cite{harrdiri}.

Metamath is a formal proof language and computer verification software developed for the purpose of formalizing mathematics in a minimalistic foundational theory \cite{metamath}. On May 12, 2016 and June 1, 2016 respectively, the author formally verified the following theorems in the Metamath formal system:
\begin{theorem} [\textmd{\thlink{dirith}}, Dirichlet's theorem]\label{thm:dirith}
$$N\in\mathbb N\land A\in\mathbb Z\land\gcd(A,N)=1\ \to\ \{p\in\mathbb P\mid N|(p-A)\}\approx\mathbb N$$
\end{theorem}
\begin{theorem} [\textmd{\thlink{pnt}}, The prime number theorem]\label{thm:pnt}
$$\left(x\in(1,\infty)\mapsto \frac{\pi(x)}{x/\log x}\right)\rightsquigarrow 1$$
\end{theorem}

The latter expression is Metamath's notation for $\lim_{x\to\infty}\frac{\pi(x)}{x/\log x}=1$, restricted to $x>1$, which is the domain of definition of the function.

These two theorems are interesting formalization targets as they both have simple statements and ``deep'' proofs, and they are also both members of the ``Formalizing 100 theorems'' list maintained by Freek Wiedijk \cite{100}, which tracks formalizations of 100 of the most famous theorems in mathematics.

Both proofs were written concurrently, over the course of about seven weeks between April 7 and June 1, 2016. This was done mostly because both theorems are in the same general subject (elementary number theory) and required similar techniques (mostly asymptotic approximation of finite sums of reals). The primary informal text used for the proof was Shapiro \cite{shapiro}, which devotes a section to Dirichlet's theorem and the whole final chapter to Selberg and Erd\H{o}s's proof of the prime number theorem.

\section{Background}\label{sec:bg}
The present work is only a broad overview of the problem and proof method. Interested readers are invited to consult the main theorems \thlink{pnt} and \thlink{dirith} at \cite{setmm}, where the exact proof is discussed in detail.

The main arithmetic functions used in the formalization are:

\begin{gather*}
\pi(x)=\left|\{p\in\mathbb P\mid p\le x\}\right|=\sum_{p\le x}1\qquad
\theta(x)=\sum_{p\le x}\log p\\
\Lambda(n)=\begin{cases}\log p&\exists p\in\mathbb P,k>0:n=p^k\\0&o.w.\end{cases}\qquad
\psi(x)=\sum_{n\le x}\Lambda(n)
\end{gather*}

Additionally, the M\"{o}bius function $\mu(n)$ is a very useful tool in sum manipulations. It is the unique multiplicative function such that $\mu(1)=1$ and $\sum_{d|n}\mu(d)=0$ for $n>1$. This yields the M\"{o}bius inversion formula: if $f(n)=\sum_{d|n}g(d)$, then $g(n)=\sum_{d|n}\mu(d)f(d)$. Since $|\mu(n)|\le 1$, this is a very powerful technique for estimating sums ``by inversion''.

The proof of Hadamard and Vall\'{e}e-Poussin relies on some deep theorems in complex analysis, such as Cauchy's theorem, which were not available at the time of this formalization, so instead we targeted the ``elementary'' proof discovered half a century later semi-independently by Erd\H{o}s and Selberg. The key step in both proofs is the Selberg symmetry formula:

\begin{theorem} [\textmd{\thlink{selberg}}, Selberg symmetry formula]\label{thm:selberg}
$$\sum_{n\le x}\Lambda(n)\log n+\sum_{uv\le x}\Lambda(u)\Lambda(v)=
2x\log x+O(x).$$
\end{theorem}

In Selberg's proof, we leverage this theorem to produce a bound on the residual $R(x)=\psi(x)-x$:

\begin{theorem} [\textmd{\thlink{pntrlog2bnd}}]\label{thm:pntrlog2bnd}
$$|R(x)|\log^2x\le2\sum_{n\le x}|R(x/n)|\log n+O(x\log x).$$
\end{theorem}

The goal is to show $\pi(x)\sim\frac x{\log x}$, but it is easily shown that $\psi(x)\sim\theta(x)\sim\pi(x)\log x$, so it is equivalent to show that $\psi(x)\sim x$, or $R(x)/x\to 0$, to establish the PNT. Given an eventual bound $|R(x)|\le ax$ and the estimation $\sum_{n\le x}\frac{\log n}n=\frac12\log^2x+O(\frac{\log x}x)$, an application of \autoref{thm:pntrlog2bnd} reproduces the original estimate $|R(x)|\le ax+o(x)$, but using improved bounds on $R(x)$ on small intervals we can improve the estimate to $|R(x)|\le (a-ca^3)x+o(x)$ for a fixed constant $c$, which produces a sequence of eventual bounds approaching zero, which proves $R(x)/x\to 0$ as desired.

In Dirichlet's theorem, the focal point is instead the Dirichlet characters $\bmod N$, which are group homomorphisms from $(\mathbb Z/N\mathbb Z)^*$ to $\mathbb C^*$, extended to $\mathbb Z/N\mathbb Z$ with value $0$ at non-units, but the general theme of estimation of sums involving $\mu,\Lambda,\log$ and the characters $\chi(n)$ is the same.

\section{Formalization}\label{sec:formal}
In keeping with Metamath's tradition of minimal complexity, we used a minimum of definition. Asymptotic estimations are reduced to the class $O(1)$ of eventually bounded functions, partial functions $\mathbb R\to\mathbb C$ such that for some $c,A$, $x\ge c$ implies $|f(x)|\le A$. An equation such as $f\in O(g)$ is rewritten as $f/g\in O(1)$ (which is correct as long as $g$ is eventually nonzero, which is always true in cases of interest), and similarly $f\in o(g)$ is rewritten as $f/g\rightsquigarrow 0$.

A few finite summation theorems take us a long way; two number-theory specific summation theorems are the following divisor sum commutations:

$$\sum_{k|n}\sum_{d|k}A(k,d)=\sum_{d|n}\sum_{m|n/d}A(dm,d)$$
$$\sum_{n\le x}\sum_{d|n}A(n,d)=\sum_{d\le x}\sum_{m\le n/d}A(dm,d)$$

A small amount of calculus was used in the proof, mostly through the following sum estimation theorem, which for example evaluates $\sum_{n\le x}\frac{\log n}n=\frac12\log^2x+O(\frac{\log x}x)$:

\begin{theorem} [\textmd{\thlink{dvfsumrlim}}]\label{thm:dvfsumrlim} If $F$ is a differentiable function with $F'=f$, and $f$ is a positive decreasing function that converges to zero, then $g(x)=\sum_{n\le x}f(n)-F(x)$ converges to some $L$ and $|g(x)-L|\le f(x)$.
\end{theorem}

\section{Comparison and Conclusion}\label{sec:conclude}
\begin{table}[ht]
\begin{center}
  \caption{Comparison of the present proof with \cite{harrdiri,avigad}. ``?'' marks an estimated or unknown value.}

  \bigskip
  
  \begin{tabular}{r|c|c|c|c|}
    & Dirichlet & PNT & Dirichlet & PNT \\
    & (author) & (author) & \cite{harrdiri} & \cite{avigad} \\ \hline
    Total time spent & 2 weeks & 5 weeks & 5 days & 12 weeks? \\ \hline
    Lines of code & 3595 & 5100 & 1183 & 19713 \\ \hline
    Compressed bytes (gzip) & 109683 & 156226 & 11762 & 97470 \\ \hline
    Informal text & 10 pp. & 37 pp. & 192 lines & 37 pp. \\ \hline
    Informal text (gzip) & 5500? & 20350? & 2524 & 20350? \\ \hline
    de Bruijn factor & 19.9? & 7.67? & 4.66 & 4.78? \\ \hline
    Verification time & 0.18 s & 0.23 s & 450 s & 1800 s? \\ \hline
  \end{tabular}
\end{center}
\end{table}

The comparison of parallel proof attempts in different systems is usually confounded by the many other factors, so these statistics should not be given undue credence. According to \cite{avigad}, Avigad's PNT project was a year-long project by four people, with the majority of the work happening during one summer, while this was a solo project over about seven work weeks. Dirichlet's theorem is 10 pages of informal text of \cite{shapiro}, and the PNT is 37 pages. Although the number of lines in the current proofs seem competitive, this is lost in the gzipped version, because the stored Metamath proof is already largely compressed, while the Isabelle and HOL scripts are plain text.

The de Bruijn factors for this work had to be estimated because the TeX source for the informal text was not available, but indications suggest that it fares poorly with comparatively large factors 19.9 and 7.67, respectively. However, when reading these statistics it is important to realize that Metamath stores {\em proofs}, not {\em proof scripts} like Isabelle and HOL. Every inference in the proof is an axiom or theorem of the system, and no proof searches are conducted by the verifier. This is reflected in the incredibly small verification time, which is normal for Metamath proofs. We do not have exact data on verification time for HOL Light, but it is believed to be on the order of minutes to hours.

These proofs are important milestones for the Metamath project. They demonstrate that even the largest of formalization projects in high level languages can also be conducted in a ``full transparency''-style system like Metamath, with entirely worked-out proofs and with all automation offloaded from the verifier to the proof generation.



\begin{thebibliography}{Abc99}
\bibitem[Dir37]{dirichlet} Dirichlet, P. G. L.: Beweis des Satzes, dass jede unbegrenzte arithmetische Progression, deren erstes Glied und Differenz ganze Zahlen ohne gemeinschaftlichen Factor sind, unendlich viele Primzahlen enth\"{a}lt. Abhand. Ak. Wiss. Berlin {\bf 48}, 313--342 (1837)

\bibitem[Sel49]{selberg} Selberg, A.: An elementary proof of the prime-number theorem. Ann. of Math. (2), Vol. 50, pp. 305--313; reprinted in Atle Selberg Collected Papers, Springer--Verlag, Berlin Heidelberg New York, 1989 {\bf 1}, 379--387 (1949)

\bibitem[Erd49]{erdos} Erd\H{o}s, P.: On a new method in elementary number theory which leads to an elementary proof of the prime number theorem. Proc. Nat. Acad. Scis. U.S.A. {\bf 35}, 374--384 (1949)

\bibitem[Avi07]{avigad} Avigad, J., Donnelly, K., Gray, D., Raff, P.: A formally verified proof of the prime number theorem. ACM Trans. Comput. Logic
{\bf 9} (1:2), 1--23 (2007)

\bibitem[Har09]{harrison} Harrison, J.: Formalizing an analytic proof of the Prime Number Theorem (dedicated to Mike Gordon on the occasion of his 60th birthday).
Journal of Automated Reasoning, 43:243--261 (2009)

\bibitem[Har10]{harrdiri} Harrison, J.: A formalized proof of Dirichlet's theorem on primes in arithmetic progression. Journal of Formalized Reasoning, [S.l.], {\bf 2} (1), 63--83 (2010)

\bibitem[Wie16]{100} Wiedijk, F.: Formalizing 100 Theorems, \url{http://www.cs.ru.nl/~freek/100/} (accessed 20 May 2016)

\bibitem[Meg07]{metamath} Megill, N.: Metamath: \textit{A Computer Language for Pure Mathematics}. Lulu Publishing, Morrisville, North Carolina (2007)

\bibitem[Sha83]{shapiro} Shapiro, H.: \textit{Introduction to the theory of numbers}. John Wiley \& Sons Inc., New York (1983)

\bibitem[Met16]{setmm} Metamath Proof Explorer, \url{http://us.metamath.org/mpegif/mmset.html} (accessed 20 May 2016)
\end{thebibliography}
\end{document}